\documentclass[lettersize,journal]{IEEEtran}
\usepackage{amsmath,amsfonts}
\usepackage{array}
\usepackage{textcomp}
\usepackage{stfloats}
\usepackage{url}
\usepackage{verbatim} 
\usepackage{graphicx}
\usepackage{cite}
\usepackage{accents}
\usepackage{amsthm}
\usepackage{xcolor}
\usepackage{tabularx, bm}
\usepackage{caption, subcaption}

\usepackage{soul} 
\usepackage{hyperref}
\hypersetup{colorlinks=true, linkcolor=blue, citecolor=blue, filecolor=blue, urlcolor=blue, pdftitle=, pdfauthor=, pdfsubject=, pdfkeywords=}
\usepackage{algorithm}
\usepackage[noend]{algpseudocode}
\usepackage{amsfonts} 

\DeclareMathOperator*{\argmin}{arg\,min}

\newtheorem{theorem}{Theorem}

\newtheorem{proposition}[theorem]{Proposition}
\newtheorem{remark}[theorem]{Remark}

\newcommand{\N}{\mathcal{N}} 
\newcommand{\E}{\mathcal{E}} 
\newcommand{\Et}{\mathcal{E}^{\tau}} 

\newcommand{\wv}{w^{\text{v}}}
\newcommand{\wc}{w^{\text{c}}}
\newcommand{\ws}{w^{\text{s}}}
\newcommand{\bc}{b^{\text{c}}}
\newcommand{\bs}{b^{\text{s}}}
\newcommand{\gs}{g^{\text{s}}}

\newcommand{\ffp}{f^{\text{p}}} 
\newcommand{\ffq}{f^{\text{q}}}
\newcommand{\UBgp}{\overline{f}^{\text{p}}} 
\newcommand{\UBgq}{\overline{f}^{\text{q}}} 
\newcommand{\LBgp}{\underline{f}^{\text{p}}} 
\newcommand{\LBgq}{\underline{f}^{\text{q}}}

\newcommand{\ddp}{d^{\text{p}}} 
\newcommand{\ddq}{d^{\text{q}}}

\newcommand{\lpp}{l^{\text{p+}}} 
\newcommand{\lpm}{l^{\text{p--}}} 
\newcommand{\lqp}{l^{\text{q+}}} 
\newcommand{\lqm}{l^{\text{q--}}} 
 
\newcommand{\Id}{I^{\text{d}}} 
\newcommand{\Ieff}{I^{\text{eff}}} 

\newcommand{\vd}{v^{\text{d}}} 
\newcommand{\vz}{v^{\text{z}}} 

\newcommand{\Nd}{\mathcal{N}^{\text{d}}} 
\newcommand{\Ns}{\mathcal{N}^{\text{s}}} 
\newcommand{\Ed}{\mathcal{E}^{\text{d}}} 

\newcommand{\vr}{v^{\text{\tiny R}}}
\newcommand{\vi}{v^{\text{\tiny I}}}

\begin{document}

\title{An ADMM-Based Approach to the Large-Scale {\color{black} Geomagnetically Induced Current} Blocker Placement Problem}

\author{
Minseok Ryu,
Arthur Barnes,
Russell Bent,
Adam Mate, and
Jose Tabarez
\thanks{M.~Ryu is with the School of Computing and Augmented Intelligence, Arizona State University, Tempe AZ (contact: mryu2@asu.edu). R.~Bent is with the Theoretical Division, Los Alamos National Laboratory, Los Alamos NM. A.~Barnes, A.~Mate, and J.~Tabarez are with the Analytics Intelligence and Technology Division, Los Alamos National Laboratory, Los Alamos NM.}

}



\maketitle

\begin{abstract} 
We propose an Alternating Direction Method of Multipliers (ADMM)-based approach for addressing the challenge of determining  placements for geomagnetically induced current (GIC) blocking devices on large-scale transmission networks. Traditional formulations of this placement problem lead to mixed-integer nonlinear programs (MINLPs) that are difficult for state-of-the-art MINLP solvers, which often struggle to exploit problem structure and fail to produce high-quality solutions on large systems. To overcome these challenges, the proposed method employs a three-block ADMM decomposition enhanced with (i) partial linearization of the GIC-related nonlinearities, supported by conditions guaranteeing exactness, and (ii) heterogeneous penalty parameters that significantly accelerate convergence to high-quality local solutions. Numerical experiments on newly constructed large-scale benchmark systems (up to 10,000 buses) demonstrate that the proposed approach outperforms existing methods by several orders of magnitude in runtime while producing solutions whose objective values closely approach the ideal all-blocker lower bound. These results highlight the effectiveness and scalability of the proposed approach for mitigation planning efforts in large transmission networks.
\end{abstract}


\section{Introduction} \label{sec:intro}

Geomagnetically induced currents (GICs) pose a significant threat to the reliability and security of bulk electric-power systems \cite{NERC2012-gmd}. Large-scale quasi-static electric fields (E-fields) induced by geomagnetic disturbances (GMDs), driven by severe space-weather activity or intentional electromagnetic pulse events, can drive substantial GICs through transmission lines and power transformer neutrals. These currents may lead to transformer core saturation, increased reactive-power demand, creation of harmonic distortion, thermal stress, and, in extreme scenarios, cascading failures and widespread outages. As highlighted in national laboratory assessments (e.g., \cite{assante2010high}), such events represent high-consequence, low-probability risks that necessitate effective mitigation strategies.

A variety of mitigation measures have been studied, including operational controls, network hardening, and equipment-level solutions (e.g., \cite{kappenman1991gic, zhu2014blocking, liang2015optimal, liang2019optimal, lu2017optimal, ryu2020algorithms, ryu2022mitigating}). Among these, transformer neutral blocking devices have received attention due to their ability to directly interrupt GIC pathways by disconnecting the neutral-to-ground connection through a capacitive interface. However, these devices are expensive and their installation produces nonlocal effects: placing a blocker at one substation can increase GIC exposure elsewhere in the grid \cite{barnes2024review}. This coupling underscores the need for optimization-driven placement strategies.

The resulting \emph{GIC Blocker Placement} (GIC-BP) problem can be formulated as a mixed-integer nonlinear program (MINLP) that integrates (i) binary placement decisions, (ii) AC optimal power flow (AC-OPF) constraints, and (iii) GIC network constraints involving absolute-value relations. This formulation is computationally challenging: the combination of combinatorial placement decisions, non-smooth absolute-value terms, and nonconvex AC equations makes the problem difficult for general-purpose MINLP solvers such as \texttt{SCIP} \cite{BestuzhevaEtal2021ZR} and \texttt{Juniper} \cite{kroger2018juniper}, which often struggle to produce high-quality solutions within practical time limits on realistically sized networks.

To address these challenges, our prior work introduced a three-block alternating direction method of multipliers (ADMM), denoted \texttt{3ADMMB} \cite{ryu2023heuristic}, which exploits the decomposable structure of the MINLP by separating blocker-placement, GIC, and OPF subproblems. While \texttt{3ADMMB} demonstrated substantially improved scalability relative to conventional MINLP solvers on small systems with up to $\approx$150 buses, two major limitations remained: (i) the nonlinear GIC subproblem becomes increasingly costly to solve with growing network size, and (ii) the use of a homogeneous penalty parameter in ADMM leads to unbalanced residuals and slow convergence.

\vspace{1mm}
\noindent\textbf{Contributions.}
This paper introduces \texttt{3ADMMB+}, a significantly enhanced ADMM-based algorithm designed to overcome the scalability limitation of existing approaches. Our main contributions are as follows:
\begin{itemize}
    \item \textbf{Linearized GIC subproblem with exactness guarantees.}  
    We apply absolute-value reformulations and McCormick relaxations to obtain a linear approximation of the GIC feasible region. We also identify conditions, related to objective weightings and ADMM penalty parameters, under which this approximation becomes \emph{exact}, meaning that the linearized and original nonlinear models yield the same solution for the GIC variables. Under these conditions, the GIC subproblem can be solved more efficiently while maintaining solution fidelity.

    \item \textbf{Heterogeneous penalty parameters and adaptive residual balancing.}   
    By imposing separate penalties for the binary and continuous consensus constraints, and updating them using a heterogeneous residual-balancing strategy, the proposed approach accounts for the differing behavior of the partial residuals. This results in considerably fewer ADMM iterations and improved numerical stability.

    \item \textbf{Superior scalability and near-ideal solution quality.}  
    Across extensive numerical experiments, \texttt{3ADMMB+} consistently outperforms \texttt{SCIP}, \texttt{Juniper}, and \texttt{3ADMMB} in both runtime and objective value. It achieves solution times that are \emph{several orders of magnitude} faster and returns high-quality plans whose objective values nearly match the {\color{black}all-blocker reference}.
    \item \textbf{New large-scale benchmark instances.}  
    We construct GIC-BP test systems derived from the \texttt{ACTIVSg500}, \texttt{ACTIVSg2k}, and \texttt{ACTIVSg10k} synthetic grids of the Texas A\&M University \cite{birchfield2016grid, birchfield2016statistical, birchfield2018creation, li2020creation}, augmented with full GIC modeling and extreme E-field scenarios. These large-scale instances represent an important contribution enabling systematic evaluation of GIC mitigation algorithms.
\end{itemize}

The remainder of this paper is organized as follows.
Section~\ref{sec:model} describes the detailed MINLP formulation for the GIC-BP problem.
Section~\ref{sec:algo} reviews the \texttt{3ADMMB} algorithm and introduces the proposed \texttt{3ADMMB+}, including linearization conditions and heterogeneous penalty design.
Section~\ref{sec:experiments} reports numerical results on large-scale test systems comparing \texttt{3ADMMB+} with state-of-the-art solvers.
Last, Section~\ref{sec:conclusion} concludes with a discussion of limitations and future research directions, including non-uniform E-field modeling, stochastic GMD uncertainty, and joint placement–operation optimization.

\section{Mathematical model} \label{sec:model} 

The GIC-BP problem is formulated as a MINLP with the following general form:
\begin{align}
\min_{z \in \mathcal{Z}} F(z) := \min_{ \substack{ (x, w) \in \mathcal{X} \\(y,w) \in \mathcal{Y}(z)  } } h(x,y,w) \label{model:MINLP}
\end{align}		  
where $\mathcal{Z} \subseteq \mathbb{B}^n$ represents the set of feasible blocker placement decisions, $\mathcal{X}$ denotes the constraint set associated with the AC optimal power flow (AC-OPF) problem, and $\mathcal{Y}(z)$ denotes the set of GIC-related constraints parameterized by the placement decision $z$. 
The variable $w$ serves as a coupling variable linking $\mathcal{X}$ and $\mathcal{Y}(z)$, while $x$ and $y$ are local variables corresponding to $\mathcal{X}$ and $\mathcal{Y}(z)$, respectively. 
The objective function $F(z)$ in~\eqref{model:MINLP} represents the optimal value of the inner minimization problem parameterized by $z$, which is a nonlinear program (NLP) that minimizes $h(x,y,w)$ subject to $(x, w) \in \mathcal{X}$ and $(y, w) \in \mathcal{Y}(z)$ for a fixed $z \in \mathcal{Z}$.

In the remainder of this section, we provide a specific formulation of the MINLP model \eqref{model:MINLP}, grounded in  \cite{barnes2024review, ryu2023heuristic}. 
Table \ref{tab-definition} summarizes the notations used throughout this work, categorized into sets, parameters, and variables for both the AC-OPF and GIC models. The network used for AC-OPF modeling is referred to as the AC network, denoted by $(\mathcal{N}, \mathcal{E})$, while the network used for GIC modeling is referred to as the DC network, denoted by $(\Nd, \Ed)$. The DC representation of the network is derived from the AC network, as detailed in \cite{ryu2023heuristic}.

\begin{table}[h!]
\footnotesize
\centering
\caption{Nomenclature}
\label{tab-definition}    
\begin{tabular}{l||l}
\hline
\multicolumn{2}{c}{\textbf{Sets and parameters}} \\ \hline
$\mathcal{G},\ \mathcal{N}, \ \mathcal{E}$ & set of generators, buses, and lines in AC network \\
$\mathcal{E}^{\tau} \subseteq \mathcal{E} $ & set of transformers \\
$\mathcal{G}_i, \mathcal{E}_i, \Et_i$ & set of generators, lines, and transformers connected to $i \in \mathcal{N}$ \\		
$\LBgp_k, \UBgp_k$ & bounds on the real power generation of $k \in \mathcal{G}$ \\
$\LBgq_k, \UBgq_k$ & bounds on the reactive power generation of $k \in \mathcal{G}$ \\
$\kappa$ & unit penalty cost for power {\color{black}im}balance at $i \in \mathcal{N}$ \\
$\ddp_i, \ddq_i$ & real and reactive power demand at $i \in \mathcal{N}$ \\
$\underline{\vr}_i, \overline{\vr}_i$, $\underline{\vi}_i, \overline{\vi}_i$ & voltage bounds at $i \in \mathcal{N}$ \\
$\gs_i, \bs_i$ & shunt conductance and susceptance at $i \in \mathcal{N}$ \\
$g_e, b_e$ & conductance, susceptance of $e \in \mathcal{E}$ \\
$\bc_e$ & line-charging susceptance of  $e \in \mathcal{E}$ \\		
$\overline{s}_e$ & apparent power limit of line $e \in \mathcal{E}$\\
$\underline{\theta}_{ij}, \overline{\theta}_{ij}$ & bounds on the phase angle difference at $e_{ij}  \in \mathcal{E}$ \\ 
$K_e$ & loss factor of transformer $e \in \mathcal{E}^{\tau}$\\
$\overline{I}_e$ & upper limit of the effective GIC on $e \in \mathcal{E}^{\tau}$\\
\hline 
$\Nd, \Ed$ & set of buses and lines in DC network \\
$\Ns \subseteq \Nd$ & set of substations \\
$\mathcal{E}^{\text{d}-}_m$ & set of incoming lines connected to $m \in \mathcal{N}^{\text{\tiny d}}$ \\ 
$\mathcal{E}^{\text{d}+}_m$ & set of outgoing lines connected to $m \in \mathcal{N}^{\text{\tiny d}}$ \\ 
$\gamma_{\ell}$ & conductance of line $\ell \in \mathcal{E}^{\text{\tiny d}}$ \\
$a_m$ & inverse of ground resistance at $m \in \mathcal{N}^{\text{\tiny d}}$  \\		
$\xi_{\ell} $ & GIC-induced voltage sources on $\ell \in \mathcal{E}^{\text{\tiny d}}$ \\ 
\hline
\multicolumn{2}{c}{\textbf{Variables}} \\ \hline		
$\vr_i$, $\vi_i$ & real and imaginary part of complex voltage at $i \in \N$ \\
$\ffp_k, \ \ffq_k$ & real and reactive power generated by $k \in \mathcal{G}$\\
$p_{ei}, \ p_{ej}$ & real power flow on $e_{ij} \in \E$ at $i \in \N$ and $j \in \N$  \\
$q_{ei}, \ q_{ej}$ & reactive power flow on $e_{ij} \in \E$ at $i \in \N$ and $j \in \N$  \\		
$ \lpp_i, \lqp_i$ & real and reactive power load shedding at $i \in \N$ \\
$ \lpm_i, \lqm_i$ & real and reactive power over-consumed at $i \in \N$ \\		
$d^{\text{qloss}}_i$ & reactive power loss due to GIC at $i \in \N$ \\		
$\Ieff_e$ & effective GIC on $e \in \Et$ \\				
\hline 
$\Id_{\ell}$ & GIC that flows on $\ell \in \Ed$ \\
$\vd_m$ & GIC-induced voltage magnitude at $m \in \Nd$  \\		
$z_m \in \{0,1\}$ & $z_m=1$ if a device is installed at $m \in \Ns$ \\
\hline
\end{tabular}
\end{table}

\subsection{AC-OPF model}
{\color{black} This subsection describes the AC-OPF component of the MINLP model. The model determines a feasible operating point of the transmission network while accounting for the additional reactive power losses induced by GICs. The decision variables include generator outputs, bus voltages, line power flows, and load imbalance variables, all of which are coupled through the standard AC power flow equations.}

\subsubsection{Operational constraints}

{\color{black}
The operational limits of the AC network are enforced through generator, voltage, and transmission constraints. Eq. \eqref{bound_gen} bounds the power generation of each generator, while Eq.~\eqref{bound_vol_mag} restricts the bus-voltage magnitudes through their rectangular components. Eq.~\eqref{thermal_limit} limits the apparent power flow on each transmission line.}
\begin{subequations}
\label{operational_constraints}  
\begin{align}
& \ffp_k \in [\LBgp_k, \UBgp_k], \quad \ffq_k \in [\LBgq_k, \UBgq_k], \ \forall k \in \mathcal{G}, \label{bound_gen} \\ 
& \vr_i  \in [\underline{\vr}_i, \overline{\vr}_i], \ \ \vi_i  \in [\underline{\vi}_i, \overline{\vi}_i] \ \forall i \in \mathcal{N},  \label{bound_vol_mag} \\
& p_{ei}^2 + q_{ei}^2 \leq  (\overline{s}_e)^2 , \ \ p_{ej}^2 + q_{ej}^2 \leq  (\overline{s}_e)^2, \ \forall e_{ij} \in \mathcal{E}.  \label{thermal_limit} 
\end{align}  
\end{subequations}

 
\subsubsection{Power flow equations}

{\color{black} The power flow equations describe the relationship between bus voltages and line power flows. Following the rectangular AC-OPF formulation, auxiliary variables $(w_i^v,w_e^c,w_e^s)$ are introduced to represent quadratic voltage terms and voltage products, resulting in the expressions (3a)–(3h). 
}
\begin{subequations}
\label{power_flow}
\begin{align} 
& p_{ei} =  g_e \wv_i - (g_e \wc_e + b_e \ws_e),  \forall e_{ij} \in \E,\\
& p_{ej} =  g_e \wv_j - (g_e \wc_e - b_e \ws_e), \forall e_{ij} \in \E, \\
& q_{ei} = -(b_e + \bc_e/2) \wv_i + (b_e \wc_e - g_e \ws_e ), \forall e_{ij} \in \E,\\
& q_{ej} = -(b_e + \bc_e/2) \wv_j + (b_e \wc_e + g_e \ws_e ), \forall e_{ij} \in \E,\\
& \wv_i = \vr_i \vr_i + \vi_i \vi_i , \ \forall i \in \N, \\
& \wc_e = \vr_i \vr_j + \vi_i \vi_j, \ \forall e_{ij} \in \E, \\
& \ws_e = \vr_j \vi_i - \vr_i \vi_j, \ \forall e_{ij} \in \E, \\
& \tan (\underline{\theta}_{ij}) \wc_e \leq \ws_e \leq \tan (\overline{\theta}_{ij}) \wc_e. \label{angle_tangent}
\end{align}
\end{subequations}

\subsubsection{Balance equations with reactive power losses by GIC}

{\color{black} The network operating point is determined by enforcing active and reactive power balance at every bus $i \in \N$:}
\begin{subequations}
\label{balance}
\begin{align}
& \sum_{e \in \mathcal{E}_i } p_{ei} = \sum_{k \in \mathcal{G}_i } \ffp_k - \ddp_i + \lpp_i - \lpm_i - \gs_i \wv_i,  \\
& \sum_{e \in \mathcal{E}_i } q_{ei} = \sum_{k \in \mathcal{G}_i } \ffq_k - \ddq_i  + \lqp_i - \lqm_i + \bs_i \wv_i - d^{\text{qloss}}_i. \label{balance_q}
\end{align}
\end{subequations}
{\color{black} 
Here, the non-negative slack variables $(\lpp_i, \lqp_i, \lpm_i, \lqm_i)$} are introduced to ensure that the AC OPF model remains feasible under all operating conditions. They also quantify the load imbalance induced by GIC through their positive values, which the optimization seeks to minimize.

{\color{black} In addition to conventional generation and demand terms, the reactive power balance in Eq.~\eqref{balance_q} explicitly accounts for the additional losses caused by GICs through the variables $d_i^{qloss}$. These losses are computed from the effective GIC values according to Eq.~\eqref{dqloss_polar}.}
\begin{subequations}
\label{dqloss_polar}
\begin{align}
& d^{\text{qloss}}_i =  \sum_{e \in \mathcal{E}^{\tau}_i} K_e \sqrt{\wv_i} \Ieff_e, \ \forall i \in \mathcal{N}, \label{dqloss_polar_1} \\
& \Ieff_e \in [0, \overline{I}_e], \ \forall e \in \Et,
\end{align}
\end{subequations}
where  
$K_e$ is a loss factor, $\overline{I}_e$ is an upper limit of the effective GIC, and
$\Ieff_e$ is the amount of effective GIC at a transformer $e \in \Et$.
{\color{black} The effective GIC variables, $\{\Ieff_e\}_{e \in \Et}$, are computed using the GIC model described in the next subsection, thereby coupling the AC-OPF model with the GIC model.}

\subsection{GIC model}
{\color{black} 
This subsection describes the GIC model that consists of three components: (i) computation of effective GICs, (ii) DC-network equations governing GIC flows, and (iii) GIC balance equations incorporating blocker placement decisions.
}

\subsubsection{Effective GIC computation}

{\color{black}
The effective GIC is used to quantify transformer exposure and serves as the coupling variable between the AC-OPF and GIC models. Since different transformer types experience GIC differently, the effective current is computed according to the corresponding winding configuration, as summarized in \eqref{eff_gic_abs}--\eqref{GSU}. 
}
\begin{subequations}
\label{eff_gic}
\begin{align}
& \Ieff_e = | \Theta_e |, \ \forall e \in \Et,  \label{eff_gic_abs}
\end{align}
where 
\begin{align}
& \Theta_e =  \big(\frac{N_h \Id_h + N_l \Id_l}{N_h} \big),   \text{ if $e$ is GWye-GWye}, \label{Gwye} \\
& \Theta_e =   \big(\frac{N_s \Id_s + N_c \Id_c}{N_s + N_c} \big),   \text{ if $e$ is GWye-GWye Auto} , \label{Auto}\\
& \Theta_e = \Id_h,   \text{ if $e$ is GWye-Delta-GSU}, \label{GSU} 
\end{align} 
where $N_h$, $N_l$, $N_s$, and $N_c$ are the number of turns in the high-side, low-side, series, and common windings, respectively.
Here, $\Id_{\ell}$ with $\ell \in \{h, l, s, c\}$ denotes the GIC variable defined on line $\ell \in \Ed$ of (see the next subsection for details).

{\color{black} To improve numerical robustness, the nonsmooth absolute-value relation is reformulated using the complementarity representation \eqref{eff_gic_smooth}--\eqref{eff_gic_smooth_3}.}
\begin{align} \label{eff_gic_smooth}  
  & s_e^+ - s_e^- = \Theta_e, \ \forall e \in \Et, \\
  & \Ieff_e = s_e^+ + s_e^-,  \ \forall e \in \Et, \label{eff_gic_smooth_2}  \\
  & s_e^- \geq 0, \ \ s_e^+ \geq 0, \ \  s_e^+ s_e^- = 0, \ \forall e \in \Et,  \label{eff_gic_smooth_3}
\end{align}
\end{subequations}
where $s_e^+$ and $s_e^-$ are slack variables.
{\color{black} We note that the bilinear equality \( s_e^{+} s_e^{-} = 0 \) in \eqref{eff_gic_smooth_3} can be replaced by the inequality \( s_e^{+} s_e^{-} \le 0 \), since both \( s_e^{+} \) and \( s_e^{-} \) are nonnegative. In our implementation, we adopt the inequality formulation to mitigate potential numerical instability arising from the bilinear equality. }
\subsubsection{GIC computation}
{\color{black} 
Once the effective GIC formulation has been established, the GIC flowing through each branch of the DC network is computed from the induced voltages and network conductances according to Eq. \eqref{gic}. The induced voltage sources represent the external geomagnetic disturbance, whereas the nodal voltages are determined by the GIC balance equations presented below.
}
The GICs are computed by:
\begin{align}
& \Id_{\ell} = \gamma_{\ell} (\vd_m - \vd_n + \xi_{\ell}), \ \forall \ell_{mn} \in \Ed, \label{gic} 
\end{align}
where $\vd$ is the GIC-induced voltage magnitude,
{\color{black} 
$\Id_{\ell}$ is the GIC flowing on the line $\ell \in \Ed$,
$\gamma_{\ell}$ is the line conductance, and
$\xi_{\ell}$ is the GMD-induced voltage source.
} 

\subsubsection{GIC balance}
Installing a blocker disconnects the transformer neutral from a substation and changes the conductance matrix. This is modeled using the GIC balance equations:
\begin{align}
& \sum_{\ell \in \mathcal{E}^{d-}_m} \Id_{\ell} - \sum_{\ell \in \mathcal{E}^{d+}_m} \Id_{\ell} = a_m \vd_m (1-z_m), \ \forall m \in \Nd, \label{gic_balance}
\end{align} 
{\color{black}
where $\mathcal{E}^{d-}_m$ and $\mathcal{E}^{d+}_m$ are sets of incoming and outgoing lines connected to $m \in \Nd$, $a_m$ is the inverse of ground resistance, and $z_m$ is the binary variable representing the placement of the blocking device, namely $z_m=1$ if a device is installed at $m \in \Nd$.
}
Eq.~\eqref{gic_balance} deduct the effect of $a_m$ if a device is installed at $m$.
For notation brevity, we set $a_m = 0$ for $m \in \Nd \setminus \Ns$ to ensure that the binary variables $z_m$ are defined only for $m \in \Ns$
{\color{black} where $\Ns$ is a set of substations.}

\subsubsection{Budget constraint}
To limit the number of blockers, we introduce the budget constraint:
\begin{subequations}
\label{budget_constraint} 
\begin{align}
& \sum_{m \in \Ns} z_m \leq V,  \\
& z_m \in \{0,1\}, \ \forall m \in \Ns.
\end{align}
\end{subequations} 
{\color{black} where $V$ is the budget for the placement.}

\subsection{Objective function} \label{sec:objective}
We adopt the objective of minimizing load imbalance:
\begin{align}
\min \kappa \sum_{i \in \mathcal{N}} (\lpp_i + \lpm_i+\lqp_i + \lqm_i)   \label{specific_obj}
\end{align}
where $\kappa$ represents the unit penalty cost {\color{black} for load imbalance}. 
Unlike conventional AC-OPF formulations, generation costs are omitted because this work focuses on high-impact, low-frequency GMD scenarios, where maintaining system operability and mitigating transformer damage are prioritized over economic dispatch~\cite{barnes2024review}.

\subsection{Connecting specific and abstract formulations}

We begin by mapping the variables of Section \ref{sec:model} in the specific formulation to those in the abstract form:
\begin{align*}
& x \gets \left[ \{ \ffp_k,  \ffq_k \}_{k \in \mathcal{G}}, \{ \vr_i, \vi_i, \lpp_i, \lpm_i, \lqp_i, \lqm_i  \}_{i \in \mathcal{N}},  \{ p_{e}, q_{e} \}_{e \in \mathcal{E}} \right],  \\
& y \gets \left[ \{ \Id_{\ell} \}_{\ell \in \Ed}, \{ \vd_m \}_{m \in \Nd}, \{s_e^+, s_e^-\}_{e \in \Et},    \right], \\
& w \gets \left[ \{\Ieff_e\}_{e \in \Et} \right], \text{ and }  z \gets \left[ \{z_m\}_{m \in \Ns} \right].
\end{align*}
Then we obtain the following mapping of constraints:
\begin{align*}
& \mathcal{X} \gets \{ \eqref{operational_constraints}-\eqref{dqloss_polar} \}, \  \mathcal{Y}(z) \gets \{ \eqref{eff_gic}-\eqref{gic_balance} \}, \text{ and }  \mathcal{Z} \gets \{ \eqref{budget_constraint} \}.
\end{align*}
Finally, the objective function is mapped as follows:
$${\color{black}h(x,y,w)=}f(x) \gets \eqref{specific_obj}.$$
 
\section{An ADMM-based approach} \label{sec:algo} 

Computing an optimal solution to the MINLP model \eqref{model:MINLP} is often computationally intractable, motivating the development of algorithms that produce high-quality local solutions within a practical time frame. 
Following this line of research, our previous work \cite{ryu2023heuristic} introduced an ADMM-based algorithm, termed the \textit{three-block alternating direction method of multipliers with binary variables} (\texttt{3ADMMB}). 
This algorithm exhibited superior performance in both computational speed and solution quality compared to traditional MINLP solvers on test instances involving power networks with fewer than 150 buses, where conventional solvers struggled to find optimal {\color{black}or even local} solutions.
{\color{black}
However, when applied to the large-scale instances considered in this work, \texttt{3ADMMB} exhibits two primary scalability limitations. First, the nonlinear GIC subproblem becomes increasingly expensive to solve as the network size grows. Second, the use of a homogeneous penalty parameter leads to unbalanced residuals and consequently slow convergence. We provide a numerical demonstration of these limitations in Appendix \ref{apx:3admmb_limitations}.
}

To address this scalability challenge, in this paper, we propose \texttt{3ADMMB+}, an enhanced version of \texttt{3ADMMB}. 
This extension incorporates several linearization techniques to lower the computational cost per iteration and employs heterogeneous penalty parameters to reduce the total number of iterations required for convergence.
As a result, \texttt{3ADMMB+} can achieve significant improvements in computation time, by orders of magnitude, while also delivering superior solution quality, attaining objective values close to the optimal case where all blockers are installed (i.e., $z_m = 1$ for all $m \in \mathcal{N}_s$ and $V = \infty$ in \eqref{budget_constraint}), but with substantially fewer blockers.

In the remainder of this section, we begin by reviewing the existing \texttt{3ADMMB} and its connection to the proposed \texttt{3ADMMB+}, followed by the acceleration schemes (i.e., linearization and heterogeneous penalty parameter techniques) designed to address scalability limitations.

\subsection{Review of \texttt{3ADMMB}}

\texttt{3ADMMB} is an ADMM-based algorithm designed to solve the following reformulation of the MINLP model \eqref{model:MINLP}, obtained by introducing {\color{black}\emph{continuous}} copies $(z^{\text{c}}, w^{\text{c}})$ of the binary variables $z$ and the linking variables $w$, respectively:
\begin{align}
\min_{\substack{z \in \mathcal{Z}, (x,w) \in \mathcal{X}}, y, z^{\text{c}}, w^{\text{c}} } \ & f(x)  \label{model:consensus}\\
\mbox{s.t.} \ &  z - z^{\text{c}} = 0, &  \qquad (\text{dual: } \lambda ) \nonumber \\
& w^{\text{c}} - w= 0, & \qquad (\text{dual: } \mu ) \nonumber\\
& (y,w^{\text{c}}) \in \mathcal{Y}(z^{\text{c}}). \nonumber
\end{align}
Applying ADMM to this reformulation yields the \texttt{3ADMMB} algorithm introduced in \cite{ryu2023heuristic}, which consists of three-block primal updates followed by a dual update.
Specifically, the primal updates involve: (i) solving a binary optimization subproblem—namely, the blocker placement subproblem—which can be efficiently solved using a greedy algorithm \cite{ryu2023heuristic}; and (ii) solving two NLP subproblems—namely corresponding to the GIC and OPF models, sequentially.
These subproblems are parameterized by the penalty parameter $\rho > 0$, which also serves as the step size for updating the dual variables $\lambda$ and $\mu$ after the primal steps.

\subsection{Proposed algorithm \texttt{3ADMMB+}}

Algorithm \ref{alg:3admmb+} presents a pseudo-code for \texttt{3ADMMB+}, which comprises three blocks of primal updates (i.e., lines 3, 4, and 5) followed by a dual update (i.e., line 6). 
In this variant, the nonlinearity in the second-block GIC subproblem is linearized to form $\mathcal{Y}^{\text{L}}(z^{\text{c}}) \supseteq \mathcal{Y}(z^{\text{c}})$ (Section \ref{sub:linear}).
{\color{black}
To promote exactness of the relaxation, we introduce
the auxiliary penalty
\begin{align}
    g(w^{\mathrm{c}})
    := \zeta  \textstyle \sum_{e\in\mathcal{E}^{\tau}}
    I_e^{\mathrm{eff}}, \label{eq:auxiliary_penalty}
\end{align}
where $\zeta>0$ is an algorithmic penalty coefficient. 
}
Additionally, the single penalty parameter $\rho$ used in \texttt{3ADMMB} is replaced with heterogeneous penalty parameters $\rho_1$ and $\rho_2$ in lines 3--6 (Section 
\ref{sub:heterogeneous}).
We note that Algorithm \ref{alg:3admmb+} becomes \texttt{3ADMMB} when $\mathcal{Y}^{\text{L}}(z^{\text{c}}) = \mathcal{Y}(z^{\text{c}})$ and $\rho = \rho_1=\rho_2$.

In the following subsections, we present the linearization techniques used to construct $\mathcal{Y}^{\text{L}}(z^{\text{c}})$ and identify the conditions under which $\mathcal{Y}^{\text{L}}(z^{\text{c}}) = \mathcal{Y}(z^{\text{c}})$. 
We also discuss the motivation for introducing heterogeneous penalty parameters, describe the corresponding termination criterion, and present residual balancing techniques used to adaptively update these penalties.

\begin{algorithm}[h!]
\caption{A pseudo-code of \texttt{3ADMMB+}}\label{alg:3admmb+}
\begin{algorithmic}[1]
\State Initialization: $t \gets 0$, $\lambda^{(t)}$, $\mu^{(t)}$, $z^{\text{c},(t)}$, and $w^{(t)}$.
\While{not terminated}
\State ($1^{\text{st}}$-block) Blocker placement subproblem:
\begin{align*}
z^{(t+1)} \gets  \argmin_{z \in \mathcal{Z}}  \ \langle \lambda^{(t)}, z \rangle + \frac{\rho_1}{2} \| z - z^{\text{c}, (t)}  \|^2.
\end{align*}
\State ($2^{\text{nd}}$-block) GIC subproblem:
\begin{align*}
& z^{\text{c}, (t+1)}, w^{\text{c}, (t+1)} \gets \hspace{-7mm}   \argmin_{ z^{\text{c}} \in \mathbb{R}^n, (y, w^{\text{c}}) \in \mathcal{Y}^{\text{L}}(z^{\text{c}})  } g(w^{\text{c}}) - \langle \lambda^{(t)}, z^{\text{c}} \rangle + \\
& \hspace{5mm}  \langle \mu^{(t)}, w^{\text{c}} \rangle + \frac{\rho_1}{2}  \| z^{(t+1)} - z^{\text{c}}   \|^2 + \frac{\rho_2}{2}\|w^{\text{c}}  - w^{(t)} \|^2  .
\end{align*}
\State ($3^{\text{rd}}$-block) OPF subproblem:
\begin{align*}
\hspace{-3mm}  & w^{(t+1)} \gets  \argmin_{(x,w) \in \mathcal{X}} \ f(x) - \langle \mu^{(t)}, w \rangle + \frac{\rho_2}{2} \| w^{c, (t+1)} - w  \|^2.    
\end{align*}
\State Dual update:
\begin{align*}
& \lambda^{(t+1)} \gets \lambda^{(t)} + \rho_1 ( z^{(t+1)} - z^{\text{c}, (t+1)} ), \\
& \mu^{(t+1)} \gets \mu^{(t)} + \rho_2 ( w^{\text{c}, (t+1)} - w^{(t+1)} ).
\end{align*}
\State Check termination criterion based on residuals.
\State Update $\rho_1$ and $\rho_2$ if needed (e.g., using the residual balancing technique).
\EndWhile
\State Output: a blocker placement decision $z^{(t+1)}$.
\end{algorithmic} 
\end{algorithm}

\subsection{Linearization of non-linear terms in the GIC model}\label{sub:linear}

{\color{black} As the network size increases, 
the nonlinearity in $\mathcal{Y}(z^{\text{c}})$ of the GIC model can be a computational bottleneck of \texttt{3ADMMB} (see Figure \ref{fig:3ADMMB-limitations} in Appendix \ref{apx:3admmb_limitations}).} 

To reduce this cost, we apply existing linearization techniques, specifically, the absolute value linearization method \cite{bertsimas1997introduction} to constraint \eqref{eff_gic_abs}, and the McCormick inequalities \cite{mccormick1976computability} to constraint \eqref{gic_balance}.
With these techniques, we construct the feasible region $\mathcal{Y}^{\text{L}}(z^{\text{c}}) \supseteq \mathcal{Y}(z^{\text{c}})$ which is a linear programming relaxation of $\mathcal{Y}(z^{\text{c}})$.
In the following, we presents several remarks, including the case where the relaxation becomes exact (i.e., when $\mathcal{Y}^{\text{L}}(z^{\text{c}}) = \mathcal{Y}(z^{\text{c}})$), which implies that the optimal solutions of the GIC model can be computed much more efficiently.


{\color{black}
\begin{proposition}[Exactness of the absolute-value reformulation]
\label{prop:absolute_exactness} 
Suppose that {\color{black}the objective function of the second-block GIC subproblem in line 4 of Algorithm \ref{alg:3admmb+}} is strictly increasing in $I_e^{\mathrm{eff}}$ over the feasible region. Then every optimal solution of the relaxation
obtained by removing $s_e^+s_e^-=0$ satisfies the complementarity constraint.
Hence, the relaxation is exact.
\end{proposition}
\begin{proof}
Suppose, by contradiction, that there exists an
optimal solution with $s_e^+>0, s_e^->0$. 
Let $\delta=\min\{s_e^+,s_e^-\}>0$.
Define $\hat s_e^+ :=s_e^+ - \delta \geq 0, \hat s_e^-:=s_e^- - \delta \geq 0$.  
Then $\hat s_e^+-\hat s_e^- = s_e^+-s_e^- = \Theta_e$, and therefore feasibility is preserved.
Moreover,
$
\hat I_e^{\mathrm{eff}}
=
\hat s_e^+ + \hat s_e^-
=
I_e^{\mathrm{eff}}-2\delta
<
I_e^{\mathrm{eff}}
$.
Since the augmented Lagrangian is strictly
increasing in $I_e^{\mathrm{eff}}$, the modified
solution yields a strictly smaller objective value,
contradicting optimality.
Therefore, the complementarity constraint must hold at optimality.
\end{proof}

\begin{remark}[Role of $\zeta$ in \eqref{eq:auxiliary_penalty}]
\label{remark_1}

The condition of Proposition~\ref{prop:absolute_exactness}
holds whenever
$
\zeta+\mu_e^{(t)}+\rho_2
\bigl(
I_e^{\mathrm{eff}}
-
I_e^{\mathrm{eff},(t)}
\bigr)
>0
$
for all $e\in\Et$.
Furthermore, if the dual multipliers are bounded,
$|\mu_e|\le M_e$,
and
$I_e^{\mathrm{eff}}\in[0,\overline I_e]$,
then a sufficient condition is
$
\zeta>
\max_{e\in\Et}
(M_e+\rho_2\overline I_e).
$
Since the dual bounds are generally unavailable in practice,
we employ a sufficiently large value of $\zeta$
and assess the quality of the relaxation empirically (see Figure~\ref{fig:constraint_violation}).
\end{remark}
}
  
\begin{remark}[McCormick inequalities] \label{remark_2}
The bilinear term in \eqref{gic_balance} can be reformulated using an auxiliary variable $\vz_m = \vd_m z^{\text{c}}_m$ for each $m \in \Nd$. The resulting nonlinearity can then be relaxed using the following McCormick inequalities:
\begin{subequations}
\label{McCormick}
\begin{align}
\vz_m &\geq \max\{ \underline{\vd}_m z^{\text{c}}_m, \ \ \vd_m - \overline{\vd}_m (1 - z^{\text{c}}_m) \}, \\
\vz_m &\leq \min\{ \overline{\vd}_m z^{\text{c}}_m, \ \ \vd_m - \underline{\vd}_m (1 - z^{\text{c}}_m) \},
\end{align}
\end{subequations}  
where $\vd_m \in [\underline{\vd}_m, \overline{\vd}_m]$ and $z^{\text{c}}_m \in [0,1]$. 
{\color{black}
The McCormick relaxation is exact whenever $z_m^{\mathrm c}\in\{0,1\}$, a standard property of McCormick envelopes involving a binary variable \cite{mccormick1976computability}. Since the quadratic penalty term $\frac{\rho_1}{2} \|z^{(t+1)}-z^{\mathrm c}\|^2$ encourages $z^{\mathrm c}$ to remain close to the binary $z^{(t+1)}$, increasing $\rho_1$ typically improves the tightness of the relaxation.
}
\end{remark}

{\color{black} These observations suggest that the proposed linearization can provide a tight approximation of the original GIC model and may recover the original nonlinear formulation when the conditions discussed above are satisfied, while significantly reducing the per-iteration computational cost as it removes the need for an NLP solver at each iteration.

We emphasize, however, that the exactness conditions presented above are sufficient but not necessary for the execution of \texttt{3ADMMB+}. Even when these conditions are not fully satisfied, the proposed approach remains a well-defined ADMM-based heuristic. In such cases, the linearized GIC model should be viewed as an approximation of the original nonlinear model, and the exactness conditions primarily affect the quality of the approximation rather than the validity of the algorithm itself. 
}

\subsection{Heterogeneous penalties} \label{sub:heterogeneous}

In addition to reducing computational cost, it is also important to decrease the total number of iterations required for convergence. 
In \texttt{3ADMMB+}, this is achieved by introducing heterogeneous penalty parameters, $\rho_1$ and $\rho_2$, in contrast to the homogeneous penalty $\rho = \rho_1 = \rho_2$ used in \texttt{3ADMMB} \cite{ryu2023heuristic}.
Note that the penalty parameter $\rho_{1}$ is enlisted to penalize primal infeasibility corresponding to $z=z^{c}$, while $\rho_{2}$ is used to penalize primal infeasibility corresponding to $I^{ac}=I^{dc}$. 

The motivation for this approach stems from the observed heterogeneity in the behavior of the partial residuals when $\rho = \rho_1 = \rho_2$ {\color{black} (see Figure \ref{fig:motivation-heterogeneous} in Appendix~\ref{apx:3admmb_limitations})}, where the partial residuals are defined as follows:
\begin{align*}
& p^{(t)}_1 := \frac{ \| z^{(t)} - z^{c,(t)} \| }{\max ( \| u^{(t)}\| , \| v^{(t)} \|)}, \ \   d^{(t)}_1 := \frac{\rho_1 \|z^{c,(t)}-z^{c,(t-1)} \|}{ \|\sigma^{(t)}\| }, \\
& p^{(t)}_2 := \frac{ \| w^{c,(t)} - w^{(t)} \| }{\max ( \| u^{(t)}\| , \| v^{(t)} \|)}, \ \  d^{(t)}_2 := \frac{\rho_2 \|w^{(t)}-w^{(t-1)} \|}{ \|\sigma^{(t)}\| },
\end{align*} 
where 
\begin{align*}
  v^{(t)} := \begin{bmatrix}
    z^{(t)} \\
    w^{c,(t)}
  \end{bmatrix}, \ \ \
   u^{(t)} := \begin{bmatrix}
    z^{c,(t)} \\
    w^{(t)}
  \end{bmatrix}, \ \ \
  \sigma^{(t)} := \begin{bmatrix}
  \lambda^{(t)}	 \\
  \mu^{(t)}
  \end{bmatrix}.  
\end{align*} 
The primal and dual residuals used in \texttt{3ADMMB} \cite{ryu2023heuristic} are defined as
\begin{align*}
p^{(t)} = \sqrt{ \left(p_1^{(t)} \right)^2 + \left(p_2^{(t)}\right)^2 }, \quad d^{(t)} = \sqrt{ \left(d_1^{(t)}\right)^2 + \left(d_2^{(t)}\right)^2 },
\end{align*}
under the homogeneous penalty (i.e., $\rho = \rho_1 = \rho_2$).
The primal residual measures constraint infeasibility, while the dual residual reflects suboptimality. 
A common termination criterion is to stop the algorithm when both residuals fall below predefined thresholds.

A practical approach to minimize the number of iterations required for termination is the residual balancing technique \cite{boyd2011distributed, wohlberg2017admm} that adaptively updates the penalty parameter $\rho$, increasing it when the primal residual is large and decreasing it when the dual residual is large. 
{\color{black}Similarly,} the heterogeneous penalties can be updated using the residual balancing technique applied separately to each pair of partial residuals:
\begin{equation}
\rho_{i}^{(t)} \leftarrow \begin{cases}\tau^{-1}\rho_{i}^{(t)} & \text{if $d_{i}^{(t)} > \beta p_{i}^{(t)}$},\\
\\
\tau\rho_{i}^{(t)} & \text{if $p_{i}^{(t)} > \beta d_{i}^{(t)}$}, 
\end{cases}
\label{res_bal}
\end{equation}
for $i=1,2$, where $\tau >1$ and $\beta > 1$.

\begin{figure*}[h!]
\centering
\setlength{\tabcolsep}{6pt}

\begin{tabular}{@{} c c c c @{}}
& \texttt{ACTIVSg500} & \texttt{ACTIVSg2k} & \texttt{ACTIVSg10k}  \\[4pt]

\raisebox{20mm}{\centering $E$=10} \hspace{-5mm} & 
\includegraphics[width=0.3\textwidth,clip]{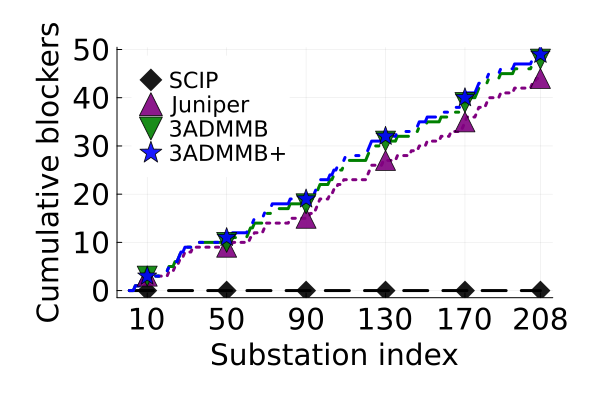} &
\includegraphics[width=0.3\textwidth,clip]{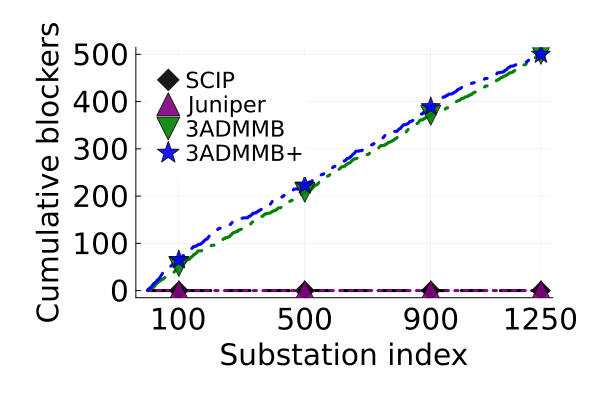} &
\includegraphics[width=0.3\textwidth,clip]{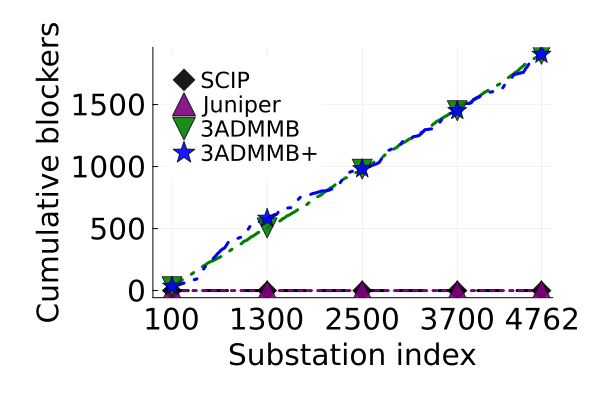} \\

\raisebox{20mm}{\centering $E$=15 } \hspace{-5mm} &
\includegraphics[width=0.3\textwidth,clip]{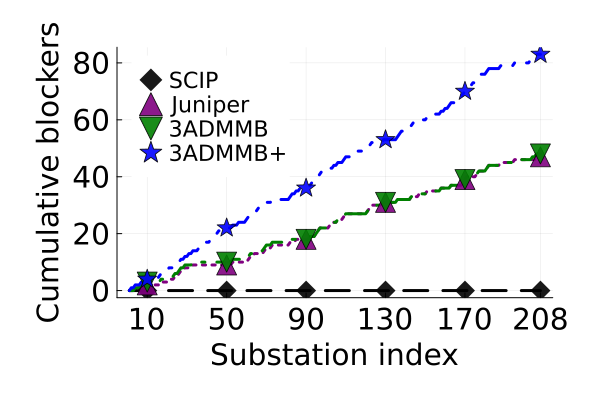} &
\includegraphics[width=0.3\textwidth,clip]{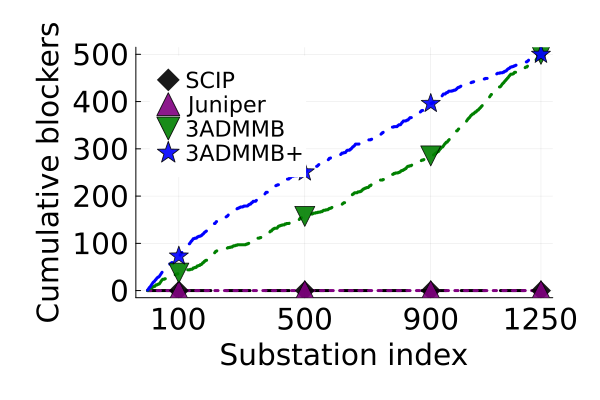} &
\includegraphics[width=0.3\textwidth,clip]{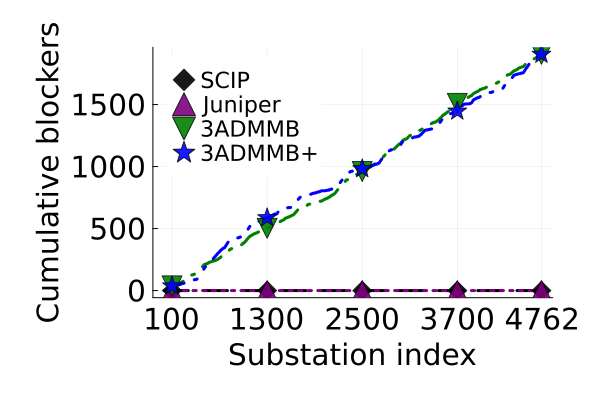} \\

\raisebox{20mm}{\centering$E$=20} \hspace{-5mm} &
\includegraphics[width=0.3\textwidth,clip]{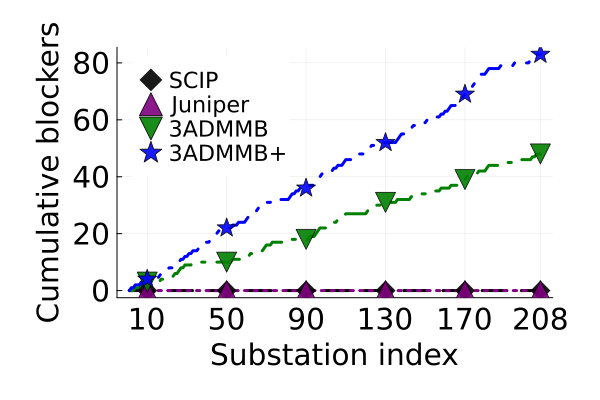} &
\includegraphics[width=0.3\textwidth,clip]{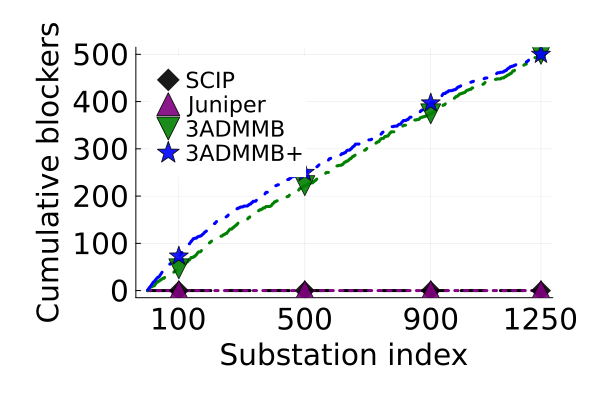} &
\includegraphics[width=0.3\textwidth,clip]{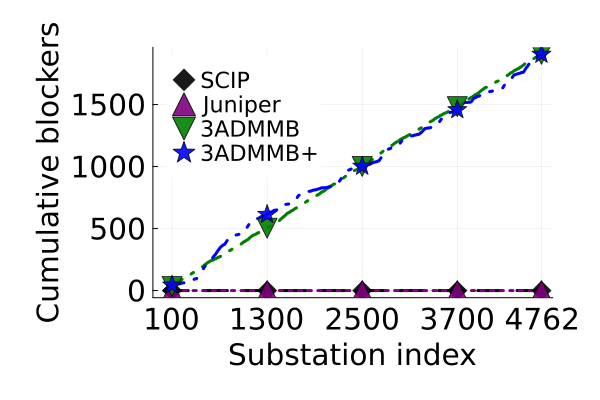} \\

\end{tabular}

\caption{Comparison of the planning decisions generated by the four approaches, $\hat{z} \in \{z^{\texttt{SCIP}}, z^{\texttt{Juniper}}, z^{\texttt{3ADMMB}}, z^{\texttt{3ADMMB+}}\}$ across different grid scales (\texttt{ACTIVSg500}, \texttt{ACTIVSg2k}, and \texttt{ACTIVSg10k})  and E-field magnitudes $E \in \{10, 15, 20\}$V/km. The x-axis in each figure corresponds to the ordered substation indices while the y-axis represents the cumulative number of substations installed with GIC blockers. For example, if the planning decision is $\hat{z} = [0, 1, 1, 0]$, meaning that substations 2 and 3 are equipped with GIC blockers out of four total substations, then the corresponding cumulative vector is $y = [0, 1, 2, 2]$.}
\label{fig:cumulative-blockers}
\end{figure*}

\begin{figure*}[h!]
\centering
\begin{subfigure}[b]{0.32\textwidth}
\centering   
\includegraphics[width=\textwidth,clip]{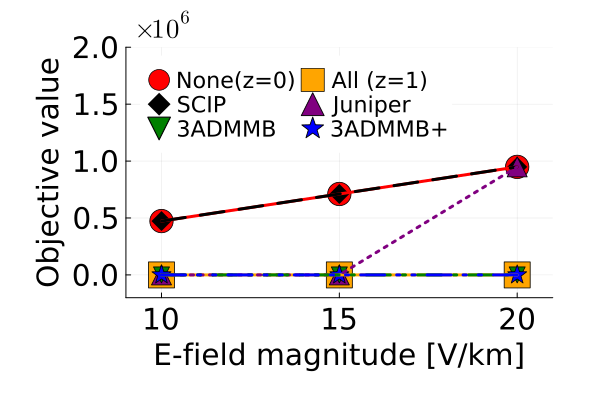} 
\includegraphics[width=\textwidth,clip]{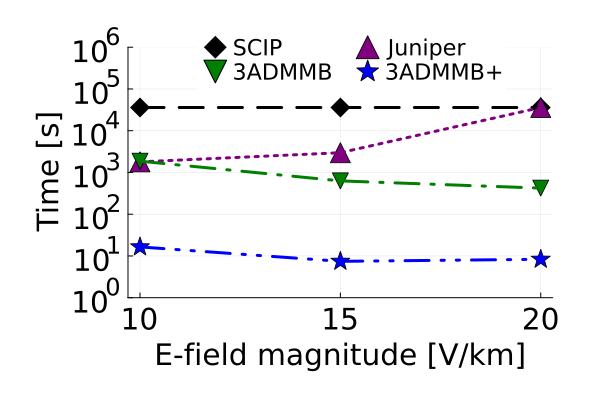} 
\caption{\texttt{ACTIVSg500}}
\label{fig:performance_comparision_500}
\end{subfigure}
\begin{subfigure}[b]{0.32\textwidth}
\centering   
\includegraphics[width=\textwidth,clip]{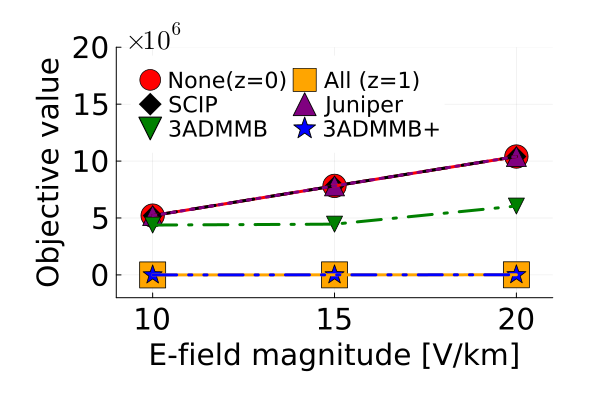} 
\includegraphics[width=\textwidth,clip]{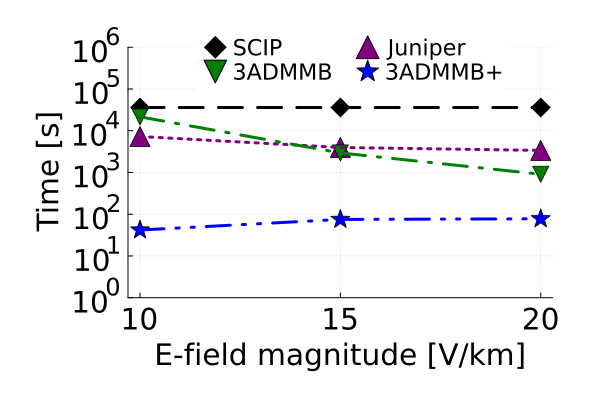} 
\caption{\texttt{ACTIVSg2k}}
\label{fig:performance_comparision_2k}
\end{subfigure}
\begin{subfigure}[b]{0.32\textwidth}
\centering   
\includegraphics[width=\textwidth,clip]{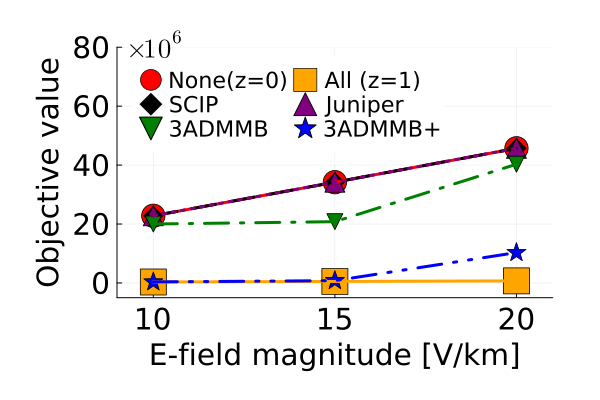} 
\includegraphics[width=\textwidth,clip]{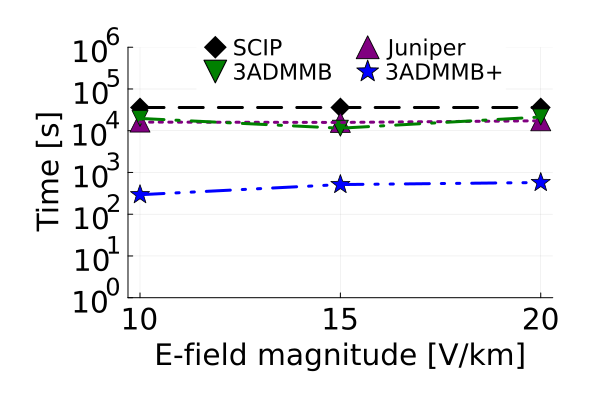} 
\caption{\texttt{ACTIVSg10k}}
\label{fig:performance_comparision_10k}
\end{subfigure} 
\caption{Performance comparison of the four approaches in terms of the objective values (i.e., $F(\hat{z})$, where $\hat{z} \in \{z^{\texttt{SCIP}}, z^{\texttt{Juniper}}, z^{\texttt{3ADMMB}}, z^{\texttt{3ADMMB+}}\}$) and computation times.}   
\label{fig:performance_comparision}
\end{figure*}

\section{Numerical experiments} \label{sec:experiments}

In this section, we numerically demonstrate the superior performance of the proposed \texttt{3ADMMB+} over existing MINLP solvers, \texttt{SCIP} and \texttt{Juniper}, and \texttt{3ADMMB} presented in \cite{ryu2023heuristic}. 
{\color{black} 
Specifically, both \texttt{SCIP} and \texttt{Juniper} are used to solve the MINLP model \eqref{model:MINLP}, where \texttt{SCIP} is a global solver that guarantees convergence to a globally optimal solution, while \texttt{Juniper} is a local solver and does not provide such a guarantee. In contrast, both \texttt{3ADMMB+} and \texttt{3ADMMB} are used to solve the reformulation \eqref{model:consensus} of the MINLP model \eqref{model:MINLP}, and likewise do not provide guarantees of global optimality.
We note that finding globally optimal solutions for the large-scale instances considered in this work within a reasonable time frame is extremely challenging. 
Therefore, we compare the solutions (i.e., planning decisions) produced by these four approaches under a fixed time limit, which can be regarded as local solutions; further details are provided in Section \ref{sec:exp_settings}.
}

Throughout the experiments, \texttt{3ADMMB+} consistently produces high-quality solutions with significantly reduced computation time, even for large-scale instances where existing approaches fail to produce meaningful planning decisions. 
{\color{black} Moreover, we observe that the solutions generated by \texttt{3ADMMB+} under the budget constraint \eqref{budget_constraint} achieve {\color{black}objective values close to those of the all-blocker reference configuration, in which all substations are equipped with GIC blocking devices}. The objective value of this reference solution is obtained by fixing all placement variables to one (i.e., $z$=1) in the MINLP model \eqref{model:MINLP} and solving the resulting continuous nonlinear programming (NLP) using IPOPT \cite{wachter2006implementation}.}

\subsection{Experimental settings}\label{sec:exp_settings}

To demonstrate the scalability of the proposed \texttt{3ADMMB+}, we generated large-scale instances based on the synthetic electrical grids \texttt{ACTIVSg500}, \texttt{ACTIVSg2k}, and \texttt{ACTIVSg10k} \cite{li2020creation}, which contain 500, 2000, and 10,000 buses, respectively.
{\color{black}
These instances are derived from the publicly available synthetic grid test cases introduced in \cite{birchfield2016grid, birchfield2016statistical, birchfield2018creation, li2020creation}. Because the original test cases do not include all the information required for GMD mitigation studies, several modifications were made. Specifically, missing data were supplemented by adding delta-wye transformers when generator step-up transformers were not explicitly modeled, specifying transformer winding configurations when omitted, and assigning geographic coordinates to substations lacking location data. 
In addition, several consistency corrections were performed, including correcting the locations of transformers modeled as connecting different substations and replacing unrealistic generator grounding resistances for generators directly connected to buses operating above 100~kV. The resulting modified test cases are publicly available through \texttt{PowerModelsGMDLib} \cite{barnes2025powermodelsgmdlib}.
}

  
To support planning decisions under the potential impacts of rare but severe geomagnetic disturbances, we consider extreme scenarios with large E-field magnitudes of $E \in \{10,15,20\}$ V/km while keeping the E-field direction fixed at 45~degrees. 
This setup is consistent with upper-bound estimates from benchmark and statistical analyses of 100-year storm events, which report E-fields in the range of $E \geq 10$ V/km depending on ground conductivity and geomagnetic latitude \cite{kelic2021literature}. 

For each scenario, we assume that at most $40\%$ of the substations can be equipped with GIC blocking devices. Accordingly, we set $V \gets \lceil 0.4  |\mathcal{N}_s| \rceil$ in~\eqref{budget_constraint}, where $|\mathcal{N}_s|$ denotes the number of substations in the grid. The resulting $(V, |\mathcal{N}_s|)$ pairs are $(83, 208)$, $(500, 1250)$, and $(1905, 4762)$ for the \texttt{ACTIVSg500}, \texttt{ACTIVSg2k}, and \texttt{ACTIVSg10k} systems, respectively.
{\color{black} To evaluate the impact of this modeling choice, Appendix~\ref{apx:budget_sensitivity} presents a sensitivity analysis over multiple budget levels, showing that the proposed approach remains effective and that the principal conclusions are not sensitive to the specific choice of a 40\% budget.}

Each instance is solved using \texttt{SCIP}, \texttt{Juniper}, \texttt{3ADMMB}, and the proposed \texttt{3ADMMB+}, producing planning decisions denoted by $z^{\texttt{SCIP}}$, $z^{\texttt{Juniper}}$, $z^{\texttt{3ADMMB}}$, and $z^{\texttt{3ADMMB+}}$, respectively, which indicate the substations where GIC blocking devices are to be installed.
The quality of each decision $\hat{z} \in \{z^{\texttt{SCIP}}, z^{\texttt{Juniper}}, z^{\texttt{3ADMMB}}, z^{\texttt{3ADMMB+}}\}$ is assessed by evaluating the objective function $F(\hat{z})$ in \eqref{model:MINLP}, which represents the penalty cost associated with the load imbalance defined in \eqref{specific_obj}, with the unit penalty costs $\kappa \gets 100$.
{\color{black}
For \texttt{3ADMMB+}, we set the auxiliary penalty coefficient $\zeta$ in \eqref{eq:auxiliary_penalty} to $10^4$.}
{\color{black} 
Given a placement decision $\hat{z}$ obtained by each method, we fix $z=\hat{z}$ in \eqref{model:MINLP} and solve the resulting NLP using IPOPT. Since the resulting NLP is nonconvex, the reported objective values correspond to locally optimal solutions obtained under this common evaluation procedure. Consequently, they provide a consistent and fair basis for comparing the planning decisions generated by different approaches.}

For the ADMM-based approaches (i.e., \texttt{3ADMMB} and \texttt{3ADMMB+}), the penalty parameters $(\rho_1, \rho_2)$ are hyperparameters that must be appropriately tuned to achieve better algorithmic performance. To mitigate this dependency, we employ the residual balancing technique \eqref{res_bal} with $\tau \gets 2$ and $\beta \gets 5$, which adaptively updates the penalty parameters based on the primal and dual residuals at each iteration.
For the initial penalty parameters, we set a relatively large value for $\rho_1$ (typically $\rho_1 \ge 10^6$) and a smaller value for $\rho_2$ (typically $\rho_2 \le 10^2$), motivated by our observation of heterogeneous behaviors in the partial residuals, as illustrated in Figure~\ref{fig:motivation-heterogeneous}.
We note that the penalty parameters $(\rho_1, \rho_2)$ may subsequently become either smaller or larger than their initial values, as they adaptively adjust according to the residual balancing technique in \eqref{res_bal}.
Under this setting, we identify the best-performing \texttt{3ADMMB} and \texttt{3ADMMB+}, respectively, for each scenario $E \in \{10, 15, 20\}$ and each grid.
For the termination criterion, we stop both algorithms when the maximum of $p^{(t)}$ and $d^{(t)}$ falls below {\color{black} the prescribed tolerance, $\epsilon=$} $5 \times 10^{-3}$.
{\color{black} Note that all NLP subproblems arising in \texttt{3ADMMB} and \texttt{3ADMMB+} are solved using IPOPT.}

Finally, we set a time limit of 10 hours for all approaches to obtain the planning decisions.
{\color{black} This limit is intended to provide a generous computational budget for all competing methods rather than to represent a practical operational requirement.}

\subsection{Performance comparison}

We compare the proposed \texttt{3ADMMB+} algorithm with \texttt{3ADMMB} \cite{ryu2023heuristic} and two MINLP solvers, \texttt{SCIP} and \texttt{Juniper}, under E-field magnitudes $E \in \{10, 15, 20\}$V/km and across grid scales \texttt{ACTIVSg500}, \texttt{ACTIVSg2k}, and \texttt{ACTIVSg10k}. Each case is tested with an placement budget limiting GIC blockers to 40\% of the substations.
Specifically, we compare the planning decisions obtained from the four approaches, $\hat{z} \in \{ z^{\texttt{SCIP}}, z^{\texttt{Juniper}}, z^{\texttt{3ADMMB}}, z^{\texttt{3ADMMB+}} \}$, along with their objective values $F(\hat{z})$ and computation times, as shown in Figures~\ref{fig:cumulative-blockers} and~\ref{fig:performance_comparision}. 

As illustrated {\color{black} in} Figure~\ref{fig:cumulative-blockers}, the MINLP solvers (\texttt{SCIP} and \texttt{Juniper}) produce nearly flat cumulative profiles, indicating that they fail to produce meaningful placement decisions within the 10-hour time limit, except for the \texttt{ACTIVSg500} cases with $E \in \{10, 15\}$ solved by \texttt{Juniper}. For these two exceptional cases, all three methods (\texttt{Juniper}, \texttt{3ADMMB}, and \texttt{3ADMMB+}) yield placement decisions achieving zero objective value, as seen on Figure~\ref{fig:performance_comparision_500}, albeit through different strategies.  

For the larger $E$=20 case on the \texttt{ACTIVSg500} grid, \texttt{Juniper} begins to return the trivial upper-bound solution $z$=0, while the two ADMM-based methods continue to produce feasible planning decisions with objective values of 59.4 and 33.3, respectively. These results indicate that \texttt{3ADMMB+} attains slightly better solution quality while requiring substantially less computation time (see Figure~\ref{fig:performance_comparision_500}).

For the larger grids, \texttt{ACTIVSg2k} and \texttt{ACTIVSg10k}, both ADMM-based approaches fully utilize the available budgets (i.e., 500 and 1905 substations, respectively), but follow slightly different deployment strategies, as evidenced by their non-overlapping cumulative curves in the second and third columns of Figure~\ref{fig:cumulative-blockers}. The superior performance of \texttt{3ADMMB+} over \texttt{3ADMMB} is corroborated by its consistently lower objective values (see the first row of Figures~\ref{fig:performance_comparision_2k} and~\ref{fig:performance_comparision_10k}), reflecting reduced GIC impacts achieved through more targeted and system-relevant blocker placements.

Figure~\ref{fig:performance_comparision} summarizes the overall performance. The proposed \texttt{3ADMMB+} offers significant improvements in both solution quality and computational efficiency. It identifies blocker placement strategies that nearly match the objective function values of the $z$=1 reference case, while reducing computation time by several orders of magnitude relative to both the MINLP solvers and \texttt{3ADMMB}. While \texttt{3ADMMB} performs better than the MINLP solvers,{\color{black}it yields objective values noticeably higher than those of the all-blocker reference and requires substantial computation time}, particularly on large-scale instances. These results demonstrate that the enhancements in \texttt{3ADMMB+} markedly improve scalability and robustness, enabling effective and timely decision-making for large-scale grids under severe geomagnetic disturbance scenarios.







\begin{figure*}[h!]
\centering
\begin{subfigure}[b]{0.24\textwidth}
\centering   
\includegraphics[width=\textwidth,clip]{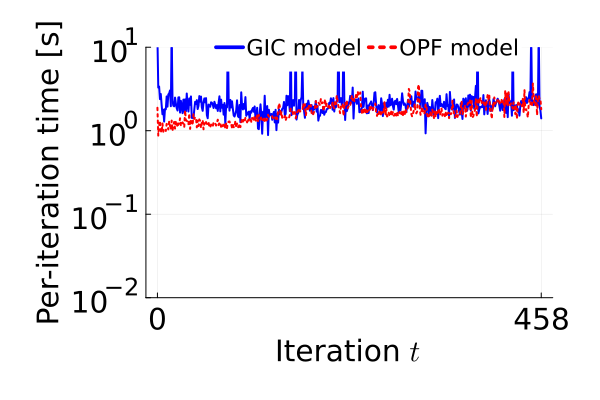}  
\includegraphics[width=\textwidth,clip]{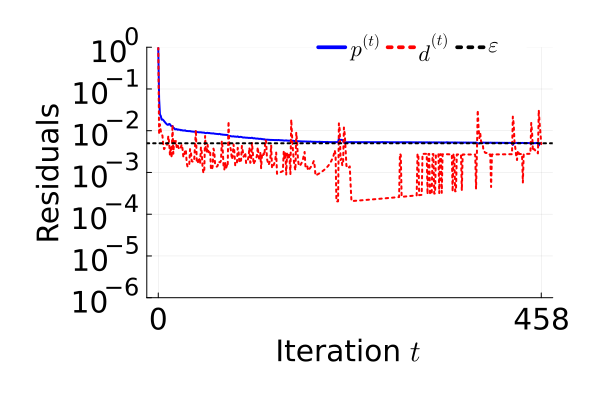}  
\caption{\texttt{3ADMMB}}
\label{fig:admm_comparison_3admmb}
\end{subfigure}  
\begin{subfigure}[b]{0.24\textwidth}
\centering   
\includegraphics[width=\textwidth,clip]{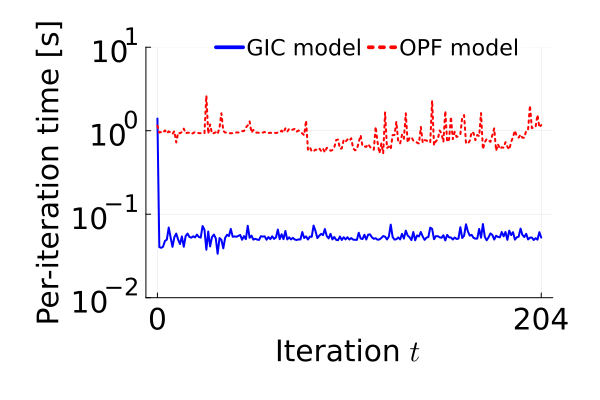}  
\includegraphics[width=\textwidth,clip]{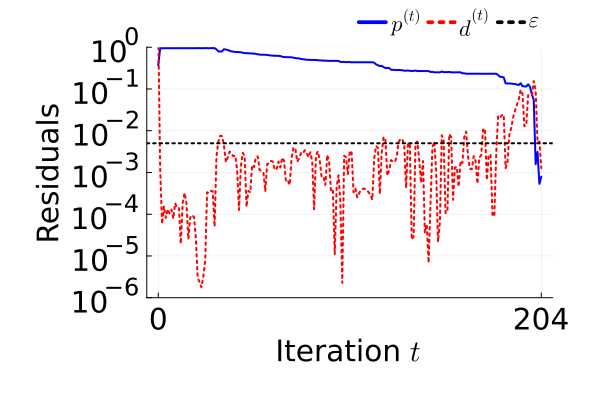}  
\caption{\texttt{3ADMMB-Linear}}
\label{fig:3admmb_mb}
\end{subfigure}  
\begin{subfigure}[b]{0.24\textwidth}
\centering   
\includegraphics[width=\textwidth,clip]{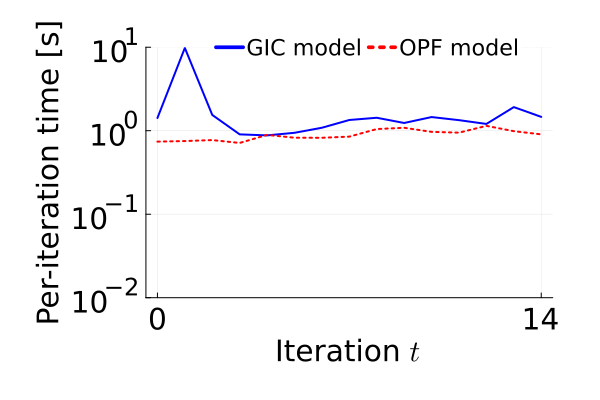}  
\includegraphics[width=\textwidth,clip]{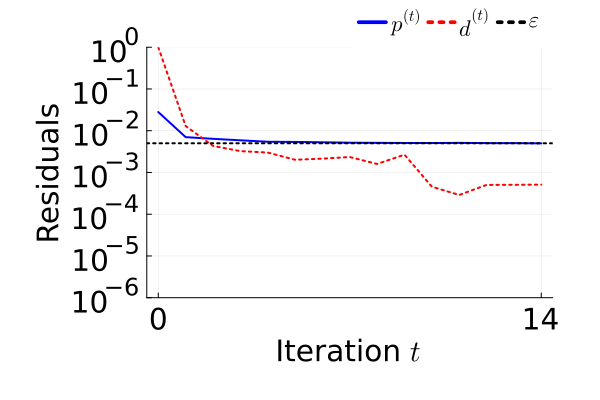}  
\caption{\texttt{3ADMMB-Hetero}}
\end{subfigure} 
\begin{subfigure}[b]{0.24\textwidth}
\centering   
\includegraphics[width=\textwidth,clip]{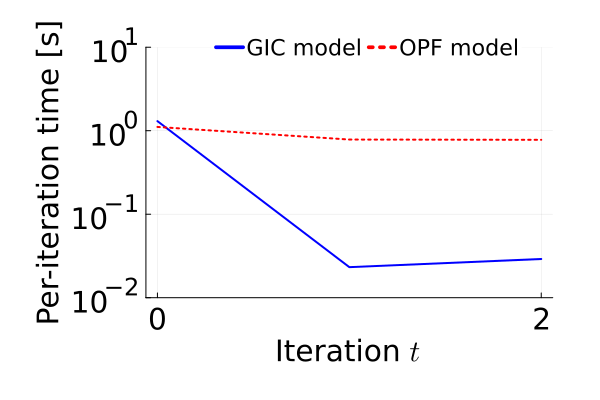}  
\includegraphics[width=\textwidth,clip]{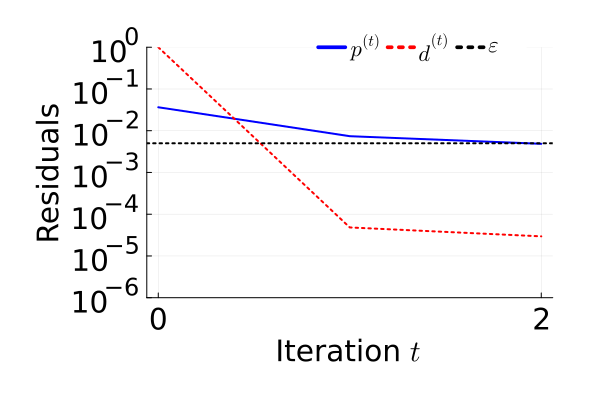}  
\caption{\texttt{3ADMMB+}}
\end{subfigure} 
\caption{Comparison of the ADMM-based approaches in terms of the per-iteration subproblem solve time {\color{black} and residual convergence}.}
\label{fig:admm_comparison}
\end{figure*}

\subsection{Acceleration schemes in \texttt{3ADMMB+}} \label{sec:acceleration}
{\color{black}
In this section, we investigate the individual contributions of the two acceleration schemes introduced in \texttt{3ADMMB+}: the linearized GIC subproblem and the heterogeneous penalty scheme. The comparison is performed on the \texttt{ACTIVSg500} grid under an $E$-field magnitude of 10 V/km, a setting in which the original \texttt{3ADMMB} produces meaningful solutions.

To isolate the effects of the two acceleration schemes, we additionally consider two intermediate variants of \texttt{3ADMMB}. The first variant, denoted by \texttt{3ADMMB-Linear}, replaces the original nonlinear GIC subproblem with its linearized counterpart while retaining the homogeneous penalty parameter used in \texttt{3ADMMB}. The second variant, denoted by \texttt{3ADMMB-Hetero}, retains the original nonlinear GIC subproblem but replaces the homogeneous penalty parameter with the proposed heterogeneous penalty scheme.

For a fair comparison, the initial penalty parameters of each algorithm are chosen as the best-performing values among the parameter settings tested, where $\rho,\rho_1\in\{10^3,10^4,10^5,10^6\}$ and $\rho_2\in\{1,10,10^2,10^3\}$. Specifically, \texttt{3ADMMB} uses $\rho=10^3$, while \texttt{3ADMMB-Linear} uses $\rho=10^5$. The heterogeneous variants employ $(\rho_1,\rho_2)=(10^3,10^2)$ for \texttt{3ADMMB-Hetero} and $(10^4,10)$ for \texttt{3ADMMB+}. For the linearized variants, we set $\zeta=10^4$, which results in only negligible complementarity violation ($\langle s^+,s^-\rangle\le10^{-6}$), as empirically verified in Appendix~\ref{apx:constraint_violation}. The larger penalty parameters used in the linearized variants promote the exactness of the linearized GIC formulation, as discussed in Section~\ref{sub:linear}. The sensitivity of \texttt{3ADMMB+} to the penalty parameters is further investigated in Appendix~\ref{apx:sensitivity}, where we show that \texttt{3ADMMB+} consistently requires substantially fewer iterations than \texttt{3ADMMB} across the tested parameter range.

Figure~\ref{fig:admm_comparison} compares the four ADMM-based algorithms in terms of per-iteration subproblem solve times and residual convergence. The first row shows that linearization reduces the computational cost of solving the GIC subproblem by more than one order of magnitude, as demonstrated by comparing \texttt{3ADMMB} with \texttt{3ADMMB-Linear}, and \texttt{3ADMMB-Hetero} with \texttt{3ADMMB+}. The OPF subproblem solve times remain essentially unchanged.

The second row compares the evolution of the primal and dual residuals. Comparing \texttt{3ADMMB} with \texttt{3ADMMB-Hetero}, and \texttt{3ADMMB-Linear} with \texttt{3ADMMB+}, shows that heterogeneous penalties substantially reduce the number of iterations required to satisfy the prescribed stopping tolerance.
{\color{black}
Linearization also reduces the iteration count even when homogeneous penalties are retained, as demonstrated by comparing \texttt{3ADMMB} with \texttt{3ADMMB-Linear}. A plausible explanation is that the enlarged feasible region of the relaxed GIC subproblem provides greater flexibility to match the consensus variables from the AC subproblem. This may contribute to the earlier reduction of the primal residual below the stopping tolerance in Figure~\ref{fig:admm_comparison}(b), while the dual residual remains small, permitting earlier termination. Thus, linearization reduces the per-iteration computational cost and may also facilitate convergence, while heterogeneous penalties provide a further reduction in iteration count.
}

Combining these two acceleration schemes, \texttt{3ADMMB+} achieves the lowest overall computation time while preserving the solution quality of the original \texttt{3ADMMB} on this test case. Both approaches produce placement decisions with zero objective value, but \texttt{3ADMMB+} requires only 5.2 seconds compared to 1863.6 seconds for \texttt{3ADMMB}, corresponding to an approximately $358\times$ speedup.

}

{\color{black}
 
\subsection{Comparison with existing heuristic approaches} \label{sec:heuristic_comparison}
We further compare the proposed \texttt{3ADMMB+} with representative heuristic approaches for GIC blocker placement. Specifically, we consider a Greedy heuristic, a simulated annealing (SA) heuristic, SA with the geographical reordering strategy proposed in~\cite{liang2019optimal} (SAR), and Greedy+SAR, in which the Greedy solution is used to initialize the SAR algorithm. 

We note that the heuristic proposed in~\cite{liang2019optimal} was originally developed for a GIC blocker placement problem without AC-OPF constraints. In contrast, our formulation explicitly incorporates AC-OPF constraints and load-shedding variables, making each candidate placement evaluation substantially more computationally expensive. Moreover, since the geographical reordering strategy was designed for a different optimization model, it is not clear \emph{a priori} whether it remains equally effective. Nevertheless, we adapted the heuristic to our formulation and evaluated its performance.
{\color{black} Further details of these heuristics are provided in Appendix~\ref{app:heuristics}.}

\begin{figure}[h!]
\centering
\includegraphics[width=0.3\textwidth,clip]{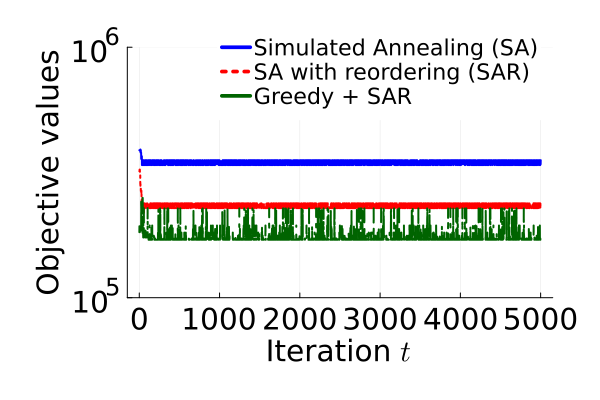}  
\caption{Comparison of heuristics on the \texttt{ACTIVSg500} instance ($E$=10V/km): Simulated annealing (SA), simulated annealing with geographical reordering (SAR), and Greedy+SAR (SAR initialized from the Greedy solution).}
\label{fig:heuristics}
\end{figure} 

Figure~\ref{fig:heuristics} compares the objective values produced by these heuristic approaches on the \texttt{ACTIVSg500} instance under an E-field magnitude of 10~V/km. As shown in the figure, the geographical reordering strategy consistently improves the original SA heuristic, and Greedy+SAR provides an additional improvement. Nevertheless, all heuristic variants remain inferior to the proposed \texttt{3ADMMB+}, which achieves a zero objective value for the same instance.
 
From the computational perspective, the heuristic approaches become increasingly expensive because every candidate blocker placement requires solving the nonlinear AC-OPF/GIC subproblem. For the \texttt{ACTIVSg500} instance, the Greedy and SA-based heuristics require approximately 208 and 6,600 seconds, respectively, whereas \texttt{3ADMMB+} computes a higher-quality solution in only approximately 10 seconds.
}

\section{Conclusion and future work} \label{sec:conclusion}

This work introduced \texttt{3ADMMB+}, an enhanced ADMM-based approach for solving the large-scale GIC-BP problem. While the original \texttt{3ADMMB} \cite{ryu2023heuristic} demonstrated the promise of ADMM-based decomposition, its high subproblem costs and slow residual convergence limited its scalability. In contrast, \texttt{3ADMMB+} achieves significant performance gains through two key enhancements: (i) linearization techniques that dramatically reduce the per-iteration computational burden and become exact under appropriate objective coefficients, and (ii) heterogeneous penalty parameters that lead to faster and more stable convergence behavior.

Across all benchmark instances, \texttt{3ADMMB+} consistently delivers solutions of substantially higher quality and at dramatically lower computational cost compared to state-of-the-art MINLP solvers (\texttt{SCIP} and \texttt{Juniper}) as well as the existing \texttt{3ADMMB}. Specifically, \texttt{3ADMMB+} improves solution times by several orders of magnitude, enabling computation on problem sizes where traditional solvers either fail to converge or exceed time limits. 
An additional contribution of this work is the construction of new large-scale GIC-BP benchmark instances, which were previously unavailable and now serve as a valuable platform for evaluating future algorithms in GIC mitigation planning.
 

The limitations of the GIC-BP problem considered in this work include the assumption of a uniform and stationary E-field on large-scale instances. 
{\color{black} Furthermore, blocker placement alone cannot completely eliminate GIC impacts under severe geomagnetic disturbances, particularly for large-scale transmission networks or high E-field magnitudes (see Appendix~\ref{apx:limitation}).}
Future extensions include incorporating non-uniform E-field models (e.g., \cite{11082653}), integrating stochastic {\color{black} or adversarial representations of E-field uncertainty \cite{ryu2020algorithms, ryu2022mitigating},} and jointly optimizing blocker placement with real-time operational controls
(e.g., line switching \cite{lu2017optimal}).

On the other hand, \texttt{3ADMMB+} is limited in that it does not guarantee global optimality for the underlying MINLP; rather, it delivers high-quality incumbent solutions. As future work, we plan to explore integrating \texttt{3ADMMB+} with local branching techniques \cite{fischetti2003local} to further improve solution quality, as well as leveraging GPU-accelerated optimization solvers (e.g., \cite{11078538}) to reduce per-iteration computational cost and enable the solution of even larger and more complex joint optimization problems under uncertainty.




\bibliographystyle{IEEEtran}
\bibliography{reference.bib}

\newpage
\section{Appendix}

{\color{black}
\subsection{Scalability limitations of \texttt{3ADMMB}}
\label{apx:3admmb_limitations}
This section provides additional numerical evidence illustrating the scalability limitations of the original \texttt{3ADMMB} algorithm that motivated the development of \texttt{3ADMMB+}.

Figure~\ref{fig:3ADMMB-limitations} reports the per-iteration solve time of the GIC and OPF subproblems in \texttt{3ADMMB}. As the network size increases, the nonlinear GIC subproblem becomes the dominant computational bottleneck. This observation motivates the linearized GIC subproblem proposed in Section~\ref{sub:linear}.

Figure~\ref{fig:motivation-heterogeneous} shows the partial primal and dual residuals associated with the two consensus constraints in \texttt{3ADMMB}. Under a homogeneous penalty parameter, the two residual pairs exhibit opposing behaviors, suggesting that different penalty parameters may be required to balance the two consensus constraints effectively. This observation motivates the heterogeneous penalty scheme proposed in Section~\ref{sub:heterogeneous}.

\begin{figure}[h!]
\centering
\includegraphics[width=0.3\textwidth,clip]{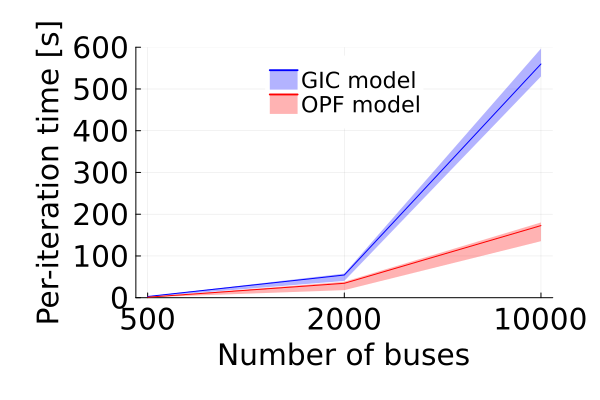}   
\caption{Average, 20th percentile, and 80th percentile of the per-iteration time for solving the GIC and OPF models—the 2nd-block and 3rd-block subproblems in \texttt{3ADMMB}, respectively.}
\label{fig:3ADMMB-limitations}
\end{figure} 

\begin{figure}[h!]
\centering
\includegraphics[width=0.3\textwidth,clip]{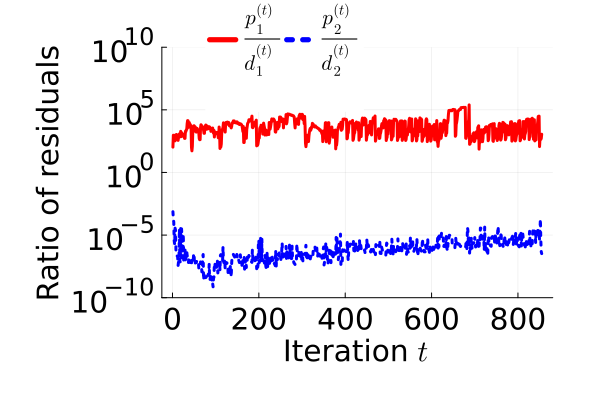}   
\caption{Opposing behaviors of the partial residuals {\color{black} in \texttt{3ADMMB}}.}
\label{fig:motivation-heterogeneous}
\end{figure}
}

{\color{black}
\subsection{Sensitivity analysis of the budget constraint}
\label{apx:budget_sensitivity}

To evaluate the impact of the budget constraint adopted in the main experiments, we conducted a sensitivity analysis by varying the budget ratio from 10\% to 90\%.

Figure~\ref{fig:budget_sensitivity} summarizes the resulting objective values produced by \texttt{3ADMMB+}.
As expected, increasing the budget monotonically decreases the objective value, reflecting the increased flexibility in placing GIC blockers. In contrast, both the runtime (approximately 5.4--5.9 s) and the number of ADMM iterations (two iterations) remain nearly unchanged across the tested budget levels. These results indicate that the computational performance of \texttt{3ADMMB+} is largely insensitive to the choice of budget and that the 40\% budget adopted in the main experiments represents a nontrivial operating point rather than a special case.

\begin{figure}[h!]
\centering
\includegraphics[width=0.3\textwidth,clip]{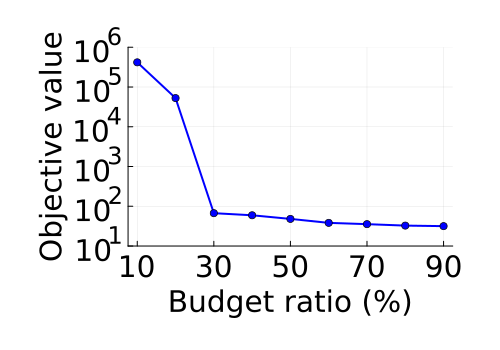}   
\caption{Sensitivity of \texttt{3ADMMB+} to the blocker budget ratio.}
\label{fig:budget_sensitivity}
\end{figure}
}

{\color{black}
\subsection{Verification of the complementarity relaxation}
\label{apx:constraint_violation}

This section empirically verifies the complementarity relaxation introduced in Section~\ref{sub:linear} (see Proposition~\ref{prop:absolute_exactness} and Remark~\ref{remark_1}). Figure~\ref{fig:constraint_violation} shows that the complementarity violation remains below $10^{-6}$ throughout the iterations of \texttt{3ADMMB+} with $\rho_1=10^4$ and $\rho_2\in\{10,100\}$. These results support the choice of $\zeta=10^4$ in our numerical experiments.

\begin{figure}[h!]
\centering
\begin{subfigure}[b]{0.3\textwidth}
\centering   
\includegraphics[width=\textwidth,clip]{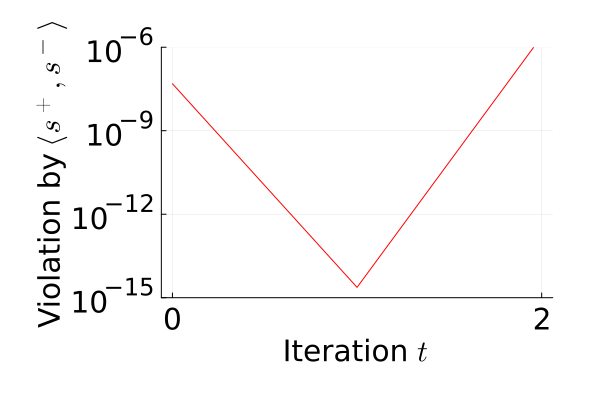}  
\caption{$\rho_2=10$} 
\end{subfigure}  
\begin{subfigure}[b]{0.3\textwidth}
\centering   
\includegraphics[width=\textwidth,clip]{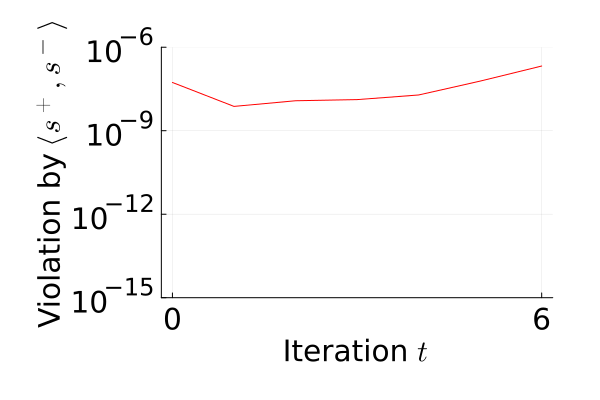}  
\caption{$\rho_2=100$} 
\end{subfigure}
\caption{Violation of $s^+_e s^-_e = 0, \ \forall e \in \Et$ in \texttt{3ADMMB+} for different value of $\rho_2 \in \{10, 100\}$, measured by $\langle s^+, s^- \rangle$.}   
\label{fig:constraint_violation}
\end{figure}

}

{\color{black} 

\begin{figure*}[h!]
\centering
\begin{subfigure}[b]{0.3\textwidth}
\centering   
\includegraphics[width=\textwidth,clip]{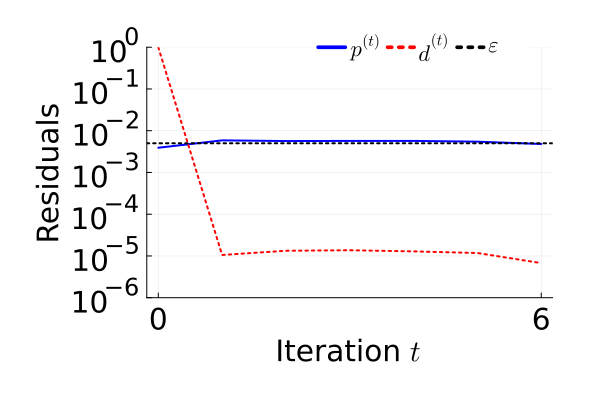}  
\caption{$\rho_2=100$}
\label{fig:penalty_parameters_2}
\end{subfigure}  
\begin{subfigure}[b]{0.3\textwidth}
\centering   
\includegraphics[width=\textwidth,clip]{./figures/8_3ADMMB_MBH_residuals_rho2_10.0.png}  
\caption{$\rho_2=10$}
\label{fig:penalty_parameters_3}
\end{subfigure}
\begin{subfigure}[b]{0.3\textwidth}
\centering   
\includegraphics[width=\textwidth,clip]{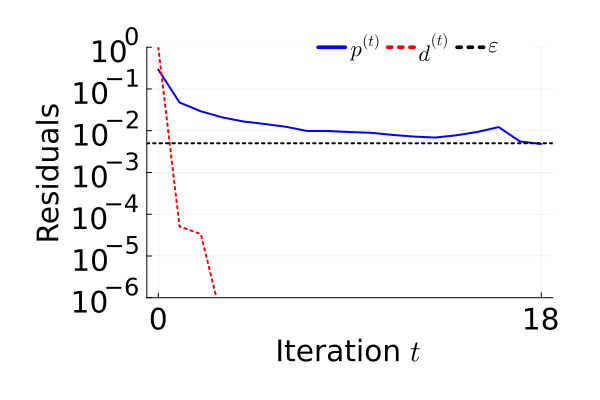}  
\caption{$\rho_2=1$}
\label{fig:penalty_parameters_4}
\end{subfigure}
\caption{Residuals of \texttt{3ADMMB+} for different values of $\rho_2 \in \{1, 10 , 100\}$.}   
\label{fig:penalty_parameters}
\end{figure*}

\subsection{Sensitivity of the heterogeneous penalty parameters}
\label{apx:sensitivity}

To investigate the sensitivity of the proposed \texttt{3ADMMB+} to the penalty parameters, we vary $\rho_2\in\{1,10,100\}$ while fixing $\rho_1=10^4$, since a sufficiently large $\rho_1$ is required to promote the exactness of the McCormick relaxation, as discussed in Remark~\ref{remark_2}.
The results are reported in Figure \ref{fig:penalty_parameters} and show that the convergence behavior of \texttt{3ADMMB+} is indeed influenced by the choice of $\rho_2$.
Nevertheless, \texttt{3ADMMB+} consistently requires substantially fewer ADMM iterations than the original \texttt{3ADMMB} (i.e., 458 iterations shown in Figure \ref{fig:admm_comparison_3admmb}) across the tested parameter range.

}
 
{\color{black} 
\subsection{Numerical illustraion of the limitation of the GIC blocker placement}
\label{apx:limitation}

{\color{black} The all-blocker configuration does not completely eliminate GIC impacts under severe geomagnetic disturbances, particularly for large-scale transmission networks or high E-field magnitudes (see Figure~\ref{fig:limitation_gic_bp}).}
This observation suggests that additional mitigation strategies, such as operational controls or transmission switching, may be required for extreme disturbance scenarios.

\begin{figure}[h!]
\centering
\includegraphics[width=0.3\textwidth,clip]{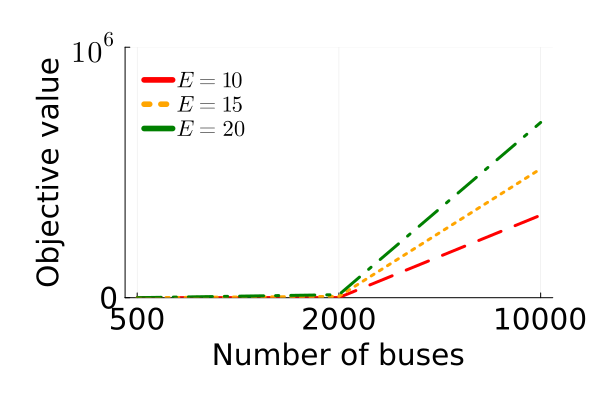}   
\caption{Limitation of the GIC-BP problem, illustrating that negative GIC impacts cannot be fully eliminated even when installing blockers at all substations, particularly as grid size and E-field magnitude increase.}
\label{fig:limitation_gic_bp}
\end{figure}
} 

{\color{black}
\subsection{Description of heuristic algorithms}
\label{app:heuristics}

We consider four heuristic algorithms for the blocker placement
problem: Greedy, simulated annealing (SA), SA with geographical
reordering (SAR), and Greedy+SAR. The SA and SAR algorithms are
adapted from \cite{liang2019optimal} to our optimization problem.

\paragraph{Greedy}
The Greedy heuristic independently evaluates each candidate
location $m \in \Ns$ by solving the optimization problem with
a single blocker installed, yielding the objective value
$F(\hat{z}_m)$, where $\hat{z}_m$ denotes the placement vector
with a blocker installed only at location $m$.
The $B$ locations with the smallest objective values are then
selected, where $B$ is the blocker budget. 

\paragraph{Simulated Annealing (SA)}
SA starts from a randomly generated feasible placement vector
$z$ satisfying the blocker budget.
At each iteration, a neighboring solution $z'$ is generated
by randomly relocating one blocker from an installed location
to an unoccupied candidate location, preserving the total
number of blockers.
Let $\Delta F = F(z') - F(z)$ denote the change in the
objective value.
The neighboring solution is accepted with probability
$
P_{\mathrm{accept}} =
\min\left\{1,\exp\left(-\frac{\Delta F}{T}\right)\right\},
$
where $T$ is the temperature parameter, which gradually
decreases during the iterations.

\paragraph{SA with geographical reordering (SAR)}
SAR uses the same acceptance criterion as SA but modifies
the neighborhood structure to encourage local moves.
We first construct the GIC DC network graph, whose vertices
and edges correspond to GMD buses and branches, respectively.
The proximity between candidate substations is measured by
the shortest-path distance, defined as the minimum number
of edges connecting them in this graph.
Starting from an initial substation, the remaining candidate
substations are sequentially ordered by selecting the nearest
unvisited substation.
Given this ordering, a neighboring placement is generated
by swapping the blocking decisions of two adjacent locations
with different blocking statuses.
This preserves the blocker budget while encouraging moves
between nearby substations.

\paragraph{Greedy+SAR}
Greedy+SAR follows the SAR procedure but uses the placement
obtained from the Greedy heuristic as its initial solution
instead of a randomly generated placement.
}

\end{document}